*Research Article*

# CONVERGENCE ANALYSIS FOR A NEW FASTER ITERATION METHOD*

**Vatan KARAKAYA[1] Yunus ATALAN*[2] Kadri DOĞAN[3]**
**Nour El Houda BOUZARA[4]**

*[1]Department of Mathematical Engineering, Yildiz Technical University, Davutpasa Campus,Esenler, 34210 Istanbul,Turkey*
*E-mail: vkkaya@yahoo.com,*
*[2,4]Department of Mathematics, Yildiz Technical University, Davutpasa Campus, Esenler, 34220 Istanbul, Turkey*
*E-mail: yunus atalan@hotmail.com, bzr.nour@gmail.com*
*[3]Department of Mathematics and Science Education, Artvin Coruh University, Artvin, Turkey*
*E-mail: dogankadri@hotmail.com*

## Abstract

In this paper, we introduce a new iteration method and show that this iteration method can be used to approximate fixed point of almost contraction mappings. Furthermore, we prove that the new iteration method is equivalent to both Mann iteration method and Picard-Mann hybrid iteration method and also converges faster than Picard-Mann hybrid iteration method for the class of almost contraction mappings. In addition to these we give a table and graphics for support this result. Finally, we prove a data dependence result for almost contraction mappings by using the new iteration method.

*Araştırma Makalesi*

## *DAHA HIZLI YENİ BİR İTERASYON METODU İÇİN YAKINSAKLIK ANALİZİ*

**Öz**

Bu makelde yeni bir iterasyon yöntemini tanımladık ve bu iterasyon yönteminin hemen hemen büzülme dönüşümlerinin sabit noktasına yaklaşımı için kullanılabilir olduğunu gösterdik. Ayrıca, yeni iterasyon yönteminin hem Mann iterasyon yöntemi hem de Picard-Mann hibrid iterasyon yöntemine denk olduğunu ve hemen hemen büzülme dönüşümleri sınıfı için Picard-Mann hibrid iterasyon yönteminden daha hızlı yakınsadığını kanıtladık. Bunlara ek olarak, bu sonucu destekleyen bir tablo ve grafik de verdik. Son olarak, yeni iterasyon yöntemini kullanarak hemen hemen büzülme dönüşümleri için bir veri bağlılığı sonucunu kanıtladık.

## 1. INTRODUCTION

Fixed point theory is one of the most useful tools of mathematics since it has different applications in many branches such as chemistry, game theory, engineering and economics etc. The fixed point of an operator $T$ is a point which is mapped to itself by an operator, that is, $x$ is a fixed point of the operator $T$ if and only if $Tx = x$.

This theory has been worked by many researches from theoretical and practical aspects. The iterative approximation of fixed point for certain classes of operators is one of the main tools in the fixed point theory. Therefore, a lot of iteration methods have been defined and studied by numerous researchers (see Berinde, 2007, Chung et al., 2012, Gürsoy et al., 2013, Karakaya et al., 2013, Mann, 1953).

Now, we give some well-known iteration methods:

A sequence $\{x_n\}_{n=0}^{\infty}$ defined by

$$\begin{cases} x_0 \in C, \\ x_{n+1} = Tx_n, \ (n \in \mathbb{N}), \end{cases} \tag{1.1}$$

is commonly known as Picard iteration method (Picard, 1890).

Mann introduced the Mann iteration method in (Mann, 1953) as follows:

$$u_{n+1} = (1-\alpha_n)u_n + \alpha_n T u_n, \ (n \in \mathbb{N}), \qquad 1.2)$$

where $0 \le \alpha_n < 1$ and $\sum_{n=0}^{\infty} \alpha_n = \infty$.

In 2007, Agarwal et al. defined the S iteration method (Agarwal et al., 2007) as

$$\begin{cases} x_{n+1} = (1-\alpha_n)Tx_n + \alpha_n T y_n \\ y_n = (1-\beta_n)x_n + \beta_n T x_n, (n \in \mathbb{N}) \end{cases} \tag{1.3}$$

The following iteration method is called Noor (Noor, 2000) iteration method:

$$\begin{cases} x_0 \in C, \\ x_{n+1} = (1-\alpha_n)x_n + \alpha_n T y_n \\ y_n = (1-\beta_n)x_n + \beta_n T z_n \\ z_n = (1-\gamma_n)x_n + \gamma_n T x_n, (n \in \mathbb{N}). \end{cases} \tag{1.4}$$

In 2013, Khan (Khan, 2013) introduced Picard-Mann hybrid iteration method as follows:

$$\begin{cases} x_0 \in C, \\ x_{n+1} = T y_n \\ y_n = (1-\alpha_n)x_n + \alpha_n T x_n, \ (n \in \mathbb{N}), \end{cases} \tag{1.5}$$

where $(\alpha_n)_{n=0}^{\infty} \subset [0,1]$.

Now, we shall introduce the following iteration method:

$$\begin{cases} x_0 \in C, \\ x_{n+1} = T[(1-\alpha_n) y_n + \alpha_n T y_n] \\ y_n = T[(1-\beta_n) x_n + \beta_n T x_n], \end{cases} \tag{1.6}$$

where $(\alpha_n)_{n=0}^{\infty}, (\beta_n)_{n=0}^{\infty} \subset [0,1]$ and $\sum_{n=0}^{\infty} \alpha_n = \infty$.

In this paper, we prove that iteration method (1.6) strongly converges to fixed point of almost contraction mappings satisfying condition (1.7). Furthermore, we show that iteration method (1.6) is equivalent to both Mann iteration method (1.2) and Picard-Mann hybrid iteration method (1.5). We also compare the rate of convergence of iteration method (1.6) and Picard-Mann hybrid iteration method (1.5). Additionally, we give a table and graphics for support this result. Finally, we give a data dependence result for almost contraction mappings satisfying condition (1.7) by employing iteration method (1.6).

In order to obtain our main results we need some lemmas and definitions given in the following:

**Lemma 1.1** (Weng, 1991). Let $\{a_n\}_{n=0}^{\infty}$ and $\{b_n\}_{n=0}^{\infty}$ be nonnegative real sequences satisfying the following condition:

$a_{n+1} \le (1-\mu_n) a_n + b_n,$

where $\mu_n \in (0,1)$, for all $n \ge n_0$, $\sum_{n=1}^{\infty} \mu_n = \infty$ and $\frac{b_n}{\mu_n} \to 0$ as $n \to \infty$. Then $\lim a_n = 0$.

**Lemma 1.2** (Soltuz, 2008). Let $\{a_n\}_{n=0}^{\infty}$ be a nonnegative real sequence and there exists an $n_0 \in \mathbb{N}$ such that for all $n \ge n_0$ satisfying the following condition:

$a_{n+1} \le (1-\mu_n) a_n + \mu_n \eta_n,$

where $\mu_n \in (0,1)$ for all $n \in \mathbb{N}, \sum_{n=1}^{\infty} \mu_n = \infty$ and $\eta_n \ge 0, \forall n \in \mathbb{N}$. Then the following inequality holds:

$0 \le \limsup_{n\to\infty} a_n \le \limsup_{n\to\infty} \eta_n.$

**Definition 1.3** (Berinde, 2003). *The self-map $T: C \to C$ is called almost contraction if there exist $\delta \in (0,1)$ and $L \ge 0$ such that*

$\|Tx - Ty\| \le \delta \|x - y\| + L \|y - Tx\|$.

**Theorem 1.4** (Berinde, 2003). *Let $(X,d)$ be a complete metric space and $T: X \to X$ be a almost contraction for which there exist $\delta \in (0,1)$ and some*

$L_1 \geq 0$ *such that*

$$\|Tx - Ty\| \leq \delta \|x - y\| + L_1 \|x - Tx\|. \tag{1.7}$$

*Then,* $T$ *has a unique fixed point.*

**Definition 1.5** (Phuengrattana and Suantai, 2013). Let $\{a_n\}_{n=1}^{\infty}$ and $\{b_n\}_{n=1}^{\infty}$ be two iterative sequences converging to the same fixed point $p_*$. We say that $\{a_n\}_{n=1}^{\infty}$ converges faster than $\{b_n\}_{n=1}^{\infty}$ to $p_*$ if

$$\lim_{n\to\infty} \frac{d(a_n, p)}{d(b_n, p)} = 0.$$

**Definition 1.6** (Soltuz, 2008). Let $T, S : C \to C$ be two operators. We say that $S$ is an approximate operator of $T$ if for all $x \in C$ and for a fixed $\varepsilon > 0$ if $\|Tx - Sx\| \leq \varepsilon.$

## 2. MAIN RESULTS

**Theorem 2.1** *Let* $C$ *be a nonempty closed convex subset of a Banach space* $X$ *and* $T : C \to C$ *be a almost contraction map satisfying condition (0.7). Let* $\{x_n\}_{n=0}^{\infty}$ *be iterative sequence generated by* (0.6) *with real sequences such that* $\{\alpha_n\}_{n=0}^{\infty}$ *and* $\{\beta_n\}_{n=0}^{\infty}$ $\in [0,1]$ *satisfying* $\sum_{k=0}^{n} \alpha_k = \infty$. *Then* $\{x_n\}_{n=0}^{\infty}$ *converges to unique fixed point* $p_*$ *of* $T$.

**Proof.** It can be seen easily from (0.7) that $p_*$ is the unique fixed point of $T$. Firstly we must show that $x_n \to p_*$ as $n \to \infty$. From (0.6) and (0.7), we have

$$\begin{aligned} \|y_n - p*\| &= \|T[(1-\beta_n)x_n + \beta_n Tx_n] - p_*\| \\ &\leq \delta \|(1-\beta_n)x_n + \beta_n Tx_n - p_*\| \\ &\leq \delta[(1-\beta_n)\|x_n - p_*\| + \beta_n \delta \|x_n - p_*\|] \\ &= \delta[1 - \beta_n(1-\delta)]\|x_n - p_*\|, \end{aligned} \tag{1.8}$$

$$\begin{aligned}\|x_{n+1}-p_*\| &= \|T[(1-\alpha_n)y_n+\alpha_n Ty_n]-p_*\| \\ &\le \delta\|(1-\alpha_n)y_n+\alpha_n Ty_n-p_*\| \\ &\le \delta[(1-\alpha_n)\|y_n-p_*\|+\alpha_n\|Ty_n-Tp_*\|] \\ &= \delta[1-\alpha_n(1-\delta)]\|y_n-p_*\| \\ &\le \delta^2[1-\alpha_n(1-\delta)][1-\beta_n(1-\delta)]\|x_n-p_*\|.\end{aligned} \tag{1.9}$$

Since $\delta\in(0,1)$ and $[1-\beta_n(1-\delta)]<1$ , we obtain

$$\|x_{n+1}-p_*\|\le[1-\alpha_n(1-\delta)]\|x_n-p_*\|. \tag{1.10}$$

By continuing this process in (0.10), we obtain the following inequalities:

$$\begin{aligned}\|x_n-p_*\| &\le [1-\alpha_{n-1}(1-\delta)]\|x_{n-1}-p_*\| \\ \|x_{n-1}-p_*\| &\le [1-\alpha_{n-2}(1-\delta)]\|x_{n-2}-p_*\| \\ &\vdots \\ \|x_1-p_*\| &\le [1-\alpha_0(1-\delta)]\|x_0-p_*\|.\end{aligned}$$

Then we have

$$\|x_{n+1}-p_*\|\le\prod_{i=0}^{n}[1-\alpha_i(1-\delta)]\|x_0-p_*\|. \tag{1.11}$$

It is well-known from the classical analysis that $1-x\le e^{-x}$ for all $x\in[0,1]$. By using this fact together with (0.11), we obtain

$$\begin{aligned}\|x_{n+1}-p_*\| &\le \prod_{i=0}^{n} e^{-(1-\delta)\alpha_i}\|x_0-p_*\| \\ &= \|x_0-p_*\|^{n+1}\left(e^{(1-\delta)\sum_{i=0}^{n}\alpha_i}\right)^{-1}.\end{aligned} \tag{1.12}$$

Taking the limit of both sides of inequality (0.12) it can be seen $x_n\to p_*$ as $n\to\infty$.

**Theorem 2.2** *Let* $X$ *be a Banach space,* $C$ *be a nonempty, closed, convex subset of* $X$ *and* $T:C\to C$ *be a almost contraction map satisfying condition* (0.7) *with fixed point* $p_*$ *. Suppose that* $\{u_n\}_{n=0}^{\infty}$ *is defined by* (0.2) *for* $u_o\in C$ *and* $\{x_n\}_{n=0}^{\infty}$ *is defined by* (0.6) *for* $x_o\in C$ *with real sequences such that* $\{\alpha_n\}_{n=0}^{\infty}$ *and*

$\{\beta_n\}_{n=0}^{\infty} \in [0,1]$. *Then the following assertions are equivalent:*

*i.* *The Mann iteration* (0.2) *converges to* $p_*$ .

*ii.* *The new iteration method* (0.6) *converges to* $p_*$ .

**Proof.** We will show that $(i) \Rightarrow (ii)$ , that is if the iteration method (0.2) converges, then the iteration method (0.6) does too. Now, by using (0.2) and (0.6) we have

$$
\begin{aligned}
\|u_{n+1} - x_{n+1}\| &= \|(1-\alpha_n)u_n + \alpha_n Tu_n - T[(1-\alpha_n)y_n + \alpha_n Ty_n]\| \\
&\le (1-\alpha_n)\|u_n - T[(1-\alpha_n)y_n + \alpha_n Ty_n]\| \\
&\quad + \alpha_n \|Tu_n - T[(1-\alpha_n)y_n + \alpha_n Ty_n]\| \\
&\le (1-\alpha_n)\{\|u_n - Tu_n\| + \|Tu_n - T[(1-\alpha_n)y_n + \alpha_n Ty_n]\|\} \\
&\quad + \alpha_n \|Tu_n - T[(1-\alpha_n)y_n + \alpha_n Ty_n]\| \\
&\le (1-\alpha_n)\|u_n - Tu_n\| + \delta\|u_n - (1-\alpha_n)y_n - \alpha_n Ty_n\| \\
&\quad + L\|u_n - Tu_n\| \qquad\qquad (1.13) \\
&\le (1-\alpha_n)\|u_n - Tu_n\| + \delta(1-\alpha_n)\|u_n - y_n\| \\
&\quad + \delta\alpha_n\|u_n - Ty_n\| + L\|u_n - Tu_n\| \\
&\le (1-\alpha_n)\|u_n - Tu_n\| + \delta(1-\alpha_n)\|u_n - y_n\| \\
&\quad + \delta\alpha_n\|u_n - Tu_n\| + \delta^2\alpha_n\|u_n - y_n\| \\
&\quad + \delta\alpha_n L\|u_n - Tu_n\| + L\|u_n - Tu_n\| \\
&= [1-\alpha_n(1-\delta) + (1+\delta\alpha_n)L]\|u_n - Tu_n\| \\
&\quad + \delta[1-\alpha_n(1-\delta)]\|u_n - y_n\|,
\end{aligned}
$$

and

$$\begin{aligned}
\|u_n - y_n\| &\le \|u_n - Tu_n\| + \|Tu_n - T[(1-\beta_n)x_n + \beta_n Tx_n]\| \\
&\le \|u_n - Tu_n\| + \delta\|u_n - (1-\beta_n)x_n - \beta_n Tx_n\| \\
&\quad + L\|u_n - Tu_n\| \\
&\le \|u_n - Tu_n\| + \delta(1-\beta_n)\|u_n - x_n\| + \delta\beta_n\|u_n - Tx_n\| \\
&\quad + L\|u_n - Tu_n\| \\
&\le \|u_n - Tu_n\| + \delta(1-\beta_n)\|u_n - x_n\| \\
&\quad + \delta\beta_n\|u_n - Tu_n\| + \delta^2\beta_n\|u_n - x_n\| \\
&\quad + \delta\beta_n L\|u_n - Tu_n\| + L\|u_n - Tu_n\| \\
&= (1+\delta\beta_n)(1+L)\|u_n - Tu_n\| \\
&\quad + \delta[1-\beta_n(1-\delta)]\|u_n - x_n\|.
\end{aligned} \tag{1.14}$$

Substituting (0.14) in (0.13), we obtain

$$\begin{aligned}
\|u_{n+1} - x_{n+1}\| &\le [1-\alpha_n(1-\delta) + (1+\delta\alpha_n)L]\|u_n - Tu_n\| \\
&\quad + \delta[1-\alpha_n(1-\delta)](1+\delta\beta_n)(1+L)\|u_n - Tu_n\| \\
&\quad + \delta^2[1-\alpha_n(1-\delta)][1-\beta_n(1-\delta)]\|u_n - x_n\| \\
&= \{[1-\alpha_n(1-\delta) + (1+\delta\alpha_n)L] \\
&\quad + \delta[1-\alpha_n(1-\delta)](1+\delta\beta_n)(1+L)\}\|u_n - Tu_n\| \\
&\quad + \delta^2[1-\alpha_n(1-\delta)][1-\beta_n(1-\delta)]\|u_n - x_n\|.
\end{aligned} \tag{1.15}$$

Since $\delta \in (0,1)$ and $[1-\beta_n(1-\delta)] < 1$, we have

$$\begin{aligned}
\|u_{n+1} - x_{n+1}\| &\le \{[1-\alpha_n(1-\delta) + (1+\delta\alpha_n)L] \\
&\quad + \delta[1-\alpha_n(1-\delta)](1+\delta\beta_n)(1+L)\}\|u_n - Tu_n\| \\
&\quad + [1-\alpha_n(1-\delta)]\|u_n - x_n\|.
\end{aligned}$$

Let

$$\begin{aligned}
\mu_n &= \alpha_n(1-\delta) \in (0,1) \\
a_n &= \|u_n - x_n\|, \\
b_n &= \{[1-\alpha_n(1-\delta) + (1+\delta\alpha_n)L] \\
&\quad + \delta[1-\alpha_n(1-\delta)](1+\delta\beta_n)(1+L)\}\|u_n - Tu_n\|.
\end{aligned}$$

Furthermore, using $Tp_* = p_*$ and $\|u_n - p_*\| \to 0$ , we have

$$\begin{aligned}\|u_n - Tu_n\| &= \|u_n - p_* + Tp_* - Tu_n\| \\ &\le \|u_n - p_*\| + \delta\|u_n - p_*\| + L\|p_* - Tp_*\| \\ &= (1+\delta)\|u_n - p_*\|.\end{aligned}$$

Then, $\|u_n - Tu_n\| \to 0$ . Because of these results, we obtain $b_n \to 0$ . By applying Lemma 1.1, we have $a_n = \|u_n - x_n\| \to 0$ as $n \to \infty$ .
Consequently;
$\|u_{n+1} - x_{n+1}\| \to 0$ as $n \to \infty$.
Now, we show that $(ii) \Rightarrow (i)$ :

$$\begin{aligned}\|x_{n+1} - u_{n+1}\| &\le \|x_{n+1} - p_*\| + \|u_{n+1} - p_*\| \\ &= \|T[(1-\alpha_n)y_n + \alpha_n Ty_n] - Tp_*\| \\ &\quad + \|(1-\alpha_n)u_n + \alpha_n Tu_n - p_*\| \\ &\le \delta\|(1-\alpha_n)y_n + \alpha_n Ty_n - p_*\| \\ &\quad + (1-\alpha_n)\|u_n - p_*\| + \alpha_n\|Tu_n - Tp_*\| \\ &\le \delta(1-\alpha_n)\|y_n - p_*\| + \delta\alpha_n\|Ty_n - Tp_*\| \\ &\quad + (1-\alpha_n)\|u_n - p_*\| + \alpha_n\delta\|u_n - p_*\| \\ &= \delta[1-\alpha_n(1-\delta)]\|y_n - p_*\| + [1-\alpha_n(1-\delta)]\|u_n - p_*\|,\end{aligned} \tag{1.16}$$

and

$$\begin{aligned}\|y_n - p_*\| &= \|T[(1-\beta_n)x_n + \beta_n Tx_n] - Tp_*\| \\ &\le \delta\|(1-\beta_n)x_n + \beta_n Tx_n - p_*\| \\ &\le \delta(1-\beta_n)\|x_n - p_*\| + \delta^2\beta_n\|x_n - p_*\| \\ &= \delta[1-\beta_n(1-\delta)]\|x_n - p_*\|.\end{aligned} \tag{1.17}$$

Substituting (0.17) in (0.16), we obtain

$$\begin{aligned}\|x_{n+1} - u_{n+1}\| &\le \delta^2[1-\alpha_n(1-\delta)][1-\beta_n(1-\delta)]\|x_n - p_*\| \\ &\quad + [1-\alpha_n(1-\delta)]\|u_n - p_*\| \\ &\le [1-\alpha_n(1-\delta)][1-\beta_n(1-\delta)]\|x_n - p_*\| \\ &\quad + [1-\alpha_n(1-\delta)]\|u_n - p_*\|.\end{aligned} \tag{1.18}$$

Denote that

$\mu_n = \alpha_n(1-\delta) \in (0,1)$

$a_n = \|u_n - p_*\|,$

$b_n = [1-\alpha_n(1-\delta)][1-\beta_n(1-\delta)]\|x_n - p_*\|.$

Thus, from Lemma 1.1, $a_n = \|u_n - p_*\| \to 0$ as $n \to \infty$. As a result of these inequalities, from (0.18),

$\|x_{n+1} - u_{n+1}\| \to 0$ as $n \to \infty.$

**Theorem 2.3** *Let* $X$ *be a Banach space,* $C$ *be a nonempty, closed, convex subset of* $X$ *and* $T: C \to C$ *be a almost contraction map satisfying condition* (0.7) *with fixed point* $p_*$ *. Suppose that* $\{u_n\}_{n=0}^{\infty}$ *is defined by* (0.5) *for* $u_o \in C$ *and* $\{x_n\}_{n=0}^{\infty}$ *is defined by* (0.6) *for* $x_o \in C$ *with real sequences such that* $\{\alpha_n\}_{n=0}^{\infty}$ *and* $\{\beta_n\}_{n=0}^{\infty}$ $\in [0,1]$ *. Then the following assertions are equivalent:*

*i.* *The Picard-Mann iteration* (0.5) *converges to* $p_*$ *.*

*ii.* *The new iteration method* (0.6) *converges to* $p_*$ *.*

**Proof.** We will show that $(i) \Rightarrow (ii)$ , that is if the iteration method (0.5) converges, then the iteration method (0.6) does too. Now, by using (0.5) and (0.6) we have

$$
\begin{aligned}
\|u_{n+1} - x_{n+1}\| &= \|Tv_n - T[(1-\alpha_n)y_n + \alpha_n Ty_n]\| \\
&\le \delta\|v_n - (1-\alpha_n)y_n - \alpha_n Ty_n\| + L\|v_n - Tv_n\| \\
&\le \delta(1-\alpha_n)\|v_n - y_n\| + \delta\alpha_n\|v_n - Ty_n\| + L\|v_n - Tv_n\| \\
&\le \delta(1-\alpha_n)\|v_n - y_n\| + (\delta\alpha_n + L)\|v_n - Tv_n\| \\
&\quad + \delta^2\alpha_n\|v_n - y_n\| + \delta\alpha_n L\|v_n - Tv_n\| \\
&= \delta[1-\alpha_n(1-\delta)]\|v_n - y_n\| + \{\delta\alpha_n(1+L) + L\}\|v_n - Tv_n\|
\end{aligned} \tag{1.19}
$$

and

$$\begin{aligned}
\|v_n - y_n\| &\le \|(1-\alpha_n)u_n + \alpha_n Tu_n - T[(1-\beta_n)x_n + \beta_n Tx_n]\| \\
&\le (1-\alpha_n)\|u_n - Tu_n\| + (1-\alpha_n)\|Tu_n - T[(1-\beta_n)x_n + \beta_n Tx_n]\| \\
&\quad +\alpha_n\|Tu_n - T[(1-\beta_n)x_n + \beta_n Tx_n]\| \\
&\le (1-\alpha_n)\|u_n - Tu_n\| + \delta\|u_n - (1-\beta_n)x_n - \beta_n Tx_n\| \\
&\quad +L\|u_n - Tu_n\| \qquad (1.20)\\
&\le (1-\alpha_n + L)\|u_n - Tu_n\| + \delta(1-\beta_n)\|u_n - x_n\| + \delta\beta_n\|u_n - Tx_n\| \\
&\le (1-\alpha_n + L)\|u_n - Tu_n\| + \delta(1-\beta_n)\|u_n - x_n\| + \delta\beta_n\|u_n - Tu_n\| \\
&\quad +\delta^2\beta_n\|u_n - x_n\| + \delta\beta_n L\|u_n - Tu_n\| \\
&= \delta[1-\beta_n(1-\delta)]\|u_n - x_n\| + \{1-\alpha_n + L + \delta\beta_n(1+L)\}\|u_n - Tu_n\|.
\end{aligned}$$

Substituting (0.20) in (0.19), we obtain

$$\begin{aligned}
\|u_{n+1} - x_{n+1}\| &\le \delta^2[1-\alpha_n(1-\delta)][1-\beta_n(1-\delta)]\|u_n - x_n\| \\
&\quad +\delta[1-\alpha_n(1-\delta)]\{1-\alpha_n + L + \delta\beta_n(1+L)\}\|u_n - Tu_n\| \\
&\quad +\{\delta\alpha_n(1+L) + L\}\|v_n - Tv_n\| \\
&\le [1-\alpha_n(1-\delta)]\|u_n - x_n\| \qquad (1.21)\\
&\quad +\delta[1-\alpha_n(1-\delta)]\{1-\alpha_n + L + \delta\beta_n(1+L)\}\|u_n - Tu_n\| \\
&\quad +\{\delta\alpha_n(1+L) + L\}\|v_n - Tv_n\|
\end{aligned}$$

Since $\delta \in (0,1)$ and $[1-\beta_n(1-\delta)] < 1$, we have

$$\begin{aligned}
\|u_{n+1} - x_{n+1}\| &\le [1-\alpha_n(1-\delta)]\|u_n - x_n\| \\
&\quad +\delta[1-\alpha_n(1-\delta)]\{1-\alpha_n + L + \delta\beta_n(1+L)\}\|u_n - Tu_n\| \\
&\quad +\{\delta\alpha_n(1+L) + L\}\|v_n - Tv_n\|.
\end{aligned}$$

Let

$$\begin{aligned}
\mu_n &= \alpha_n(1-\delta) \in (0,1) \\
a_n &= \|u_n - x_n\|, \\
b_n &= \delta[1-\alpha_n(1-\delta)]\{1-\alpha_n + L + \delta\beta_n(1+L)\}\|u_n - Tu_n\| \\
&\quad +\{\delta\alpha_n(1+L) + L\}\|v_n - Tv_n\|.
\end{aligned}$$

Furthermore, using $Tp_* = p_*$ and $\|u_n - p_*\| \to 0$, we have

$$\begin{aligned}
\|u_n - Tu_n\| &= \|u_n - p_* + Tp_* - Tu_n\| \\
&\le \|u_n - p_*\| + \delta\|u_n - p_*\| + L\|p_* - Tp_*\| \\
&= (1+\delta)\|u_n - p_*\|.
\end{aligned}$$

and

$$
\begin{aligned}
\|v_n - Tv_n\| &\le \|v_n - p_*\| + \delta\|v_n - p_*\| + L\|p_* - Tp_*\| \\
&= (1+\delta)\|v_n - p_*\| \\
&\le (1+\delta)(1-\alpha_n)\|u_n - p_*\| + (1+\delta)\alpha_n\delta\|u_n - p_*\| \\
&= (1+\delta)\left[1-\alpha_n(1-\delta)\right]\|u_n - p_*\|
\end{aligned}
$$

Then, $\|u_n - Tu_n\| \to 0$ and $\|v_n - Tv_n\| \to 0$. Because of these results, we obtain $b_n \to 0$ . By applying Lemma 1.1, we have $a_n = \|u_n - x_n\| \to 0$ as $n \to \infty$ . Consequently

$\|u_{n+1} - x_{n+1}\| \to 0$ as $n \to \infty.$

Now, we show that $(ii) \Rightarrow (i)$ :

$$
\begin{aligned}
\|x_{n+1} - u_{n+1}\| &\le \|T[(1-\alpha_n)y_n + \alpha_n Ty_n] - Tv_n\| \\
&\le \delta\|(1-\alpha_n)y_n + \alpha_n Ty_n - v_n\| \\
&\quad + L\|(1-\alpha_n)y_n + \alpha_n Ty_n - T[(1-\alpha_n)y_n + \alpha_n Ty_n]\| \\
&\le \delta(1-\alpha_n)\|y_n - u_n\| + \delta\alpha_n\|Ty_n - Tu_n\| \\
&\quad + L(1-\alpha_n)\|y_n - T[(1-\alpha_n)y_n + \alpha_n Ty_n]\| \\
&\quad + L\alpha_n\|Ty_n - T[(1-\alpha_n)y_n + \alpha_n Ty_n]\| \\
&\le \delta(1-\alpha_n)\|y_n - u_n\| + \delta^2\alpha_n\|y_n - u_n\| + \delta\alpha_n L\|y_n - Ty_n\| \\
&\quad + L(1-\alpha_n)\|y_n - Ty_n\| \\
&\quad + L(1-\alpha_n)\|Ty_n - T[(1-\alpha_n)y_n + \alpha_n Ty_n]\| \\
&\quad + L\alpha_n\|Ty_n - T[(1-\alpha_n)y_n + \alpha_n Ty_n]\| \\
&= \delta\left[1-\alpha_n(1-\delta)\right]\|y_n - u_n\| \\
&\quad + \left\{L(2\delta\alpha_n + L + 1 - \alpha_n\right\}\|y_n - Ty_n\|
\end{aligned}
\tag{1.22}
$$

and

$$
\begin{aligned}
\|y_n - u_n\| &\le \|T[(1-\beta_n)x_n + \beta_n Tx_n] - Tx_n\| + \|Tx_n - x_n\| \\
&\quad + \|x_n - u_n\| \\
&\le \delta\|(1-\beta_n)x_n + \beta_n Tx_n - x_n\| + L\|x_n - Tx_n\| \\
&\quad + \|x_n - Tx_n\| + \|x_n - u_n\| \\
&= \|x_n - u_n\| + \left\{\delta\beta_n + L + 1\right\}\|x_n - Tx_n\|.
\end{aligned}
\tag{1.23}
$$

Substituting (0.23) in (0.22), we obtain

$$\begin{aligned}\|x_{n+1}-u_{n+1}\| &\le \delta[1-\alpha_n(1-\delta)]\|x_n-u_n\| \\ &+\delta[1-\alpha_n(1-\delta)](\delta\beta_n+L+1)\|x_n-Tx_n\| \\ &+\{L(2\delta\alpha_n+L+1-\alpha_n\}\|y_n-Ty_n\|\end{aligned}$$

Denote that

$$\begin{aligned}\mu_n &= \alpha_n(1-\delta)\in(0,1) \\ a_n &= \|x_n-u_n\|, \\ b_n &= \delta[1-\alpha_n(1-\delta)](\delta\beta_n+L+1)\|x_n-Tx_n\| \\ &+\{L(2\delta\alpha_n+L+1-\alpha_n\}\|y_n-Ty_n\|.\end{aligned}$$

Furthermore, using $Tp_* = p_*$ and $\|u_n - p_*\| \to 0$ , we have

$$\begin{aligned}\|x_n-Tx_n\| &\le \|x_n-p_*\|+\delta\|x_n-p_*\|+L\|p_*-Tp_*\| \\ &= (1+\delta)\|x_n-p_*\|.\end{aligned}$$

and

$$\begin{aligned}\|y_n-Ty_n\| &\le \|y_n-p_*\|+\delta\|y_n-p_*\|+L\|p_*-Tp_*\| \\ &= (1+\delta)\|y_n-p_*\| \\ &\le (1+\delta)\delta\|(1-\beta_n)x_n+\beta_nTx_n-p_*\|+L\|p_*-Tp_*\| \\ &\le (1+\delta)\delta(1-\beta_n)\|x_n-p_*\|+(1+\delta)\delta^2\beta_n\|x_n-p_*\| \\ &= (1+\delta)\delta\left[1-\beta_n(1-\delta)\right]\|x_n-p_*\|.\end{aligned}$$

Then, $\|x_n - Tx_n\| \to 0$ and $\|y_n - Ty_n\| \to 0$ . Because of these results, we obtain $b_n \to 0$ . By applying Lemma 1.1, we have $a_n = \|u_n - x_n\| \to 0$ as $n \to \infty$ .

Consequently

$\|x_{n+1} - u_{n+1}\| \to 0$ as $n \to \infty$.

As a consequence of Theorem 2.2 and Theorem 2.3, we can give the following corollary:

**Corollary 2.4** *Let $X$ be a Banach space, $C$ be a nonempty, closed, convex subset of $X$ and $T : C \to C$ be a mapping satisfying condition (0.7) with fixed point $p_*$. If the initial point is the same for all iterations, then the following assertions are equivalent:*

*(i) the Picard iteration (0.1) converges to $p_*$ ,*

*(ii) the Mann iteration (0.2) converges to $p_*$ ,*

*(iii) the S iteration (0.3) converges to $p_*$ ,*

*(iv)the Noor iteration* (0.4) *converges to* $p_*$ ,

*(v) the Picard-Mann iteration* (0.5) *converges to* $p_*$ .

*(vi) the new iteration* (0.6) *converges to* $p_*$.

**Theorem 2.5** *Let* $X$ *be a Banach spaces, and* $C$ *be a closed, convex subset of* $X$ *, and* $T:C\to C$ *be a almost contraction map satisfying condition* (0.7) *with fixed point* $p_*$ *. Let* $\{\alpha_n\}$ *and* $\{\beta_n\}$ *be real sequences such that* $0<\alpha_n,\beta_n<1$ *for all* $n\in\mathbb{N}$ *. For given* $x_0=u_0\in C$ *, consider the iterative sequences* $\{u_n\}_{n=0}^{\infty}$ *and* $\{x_n\}_{n=0}^{\infty}$ *defined by* (0.5) *and* (0.6), *respectively. Then* $\{x_n\}_{n=0}^{\infty}$ *converges to* $p_*$ *faster than* $\{u_n\}_{n=0}^{\infty}$ *does.*

**Proof.** Let

$$\begin{aligned}\|x_{n+1}-p_*\| &= \|T[(1-\alpha_n)y_n+\alpha_n Ty_n]-p_*\| \\ &\le \delta[1-\alpha_n(1-\delta)]\|y_n-p_*\|\end{aligned},$$

by using (0.6), we get

$$\|x_{n+1}-p_*\|\le\delta^2[1-\alpha_n(1-\delta)][1-\beta_n(1-\delta)]\|x_n-p_*\|.$$

Repeating this process $n$-times we get

$$\|x_{n+1}-p_*\|\le\delta^{2(n+1)}\prod_{k=0}^{n}[1-\alpha_k(1-\delta)][1-\beta_k(1-\delta)]\|x_0-p_*\|.$$

In further it is easy to see that

$$\|u_{n+1}-p_*\|\le\delta^{n+1}\prod_{k=0}^{n}[1-\alpha_k(1-\delta)]\|u_0-p_*\|,$$

$$\begin{aligned}\lim_{n\to\infty}\frac{\|x_{n+1}-p_*\|}{\|u_{n+1}-p_*\|} &= \lim_{n\to\infty}\frac{\delta^{2(n+1)}\prod_{k=0}^{n}[1-\alpha_k(1-\delta)][1-\beta_k(1-\delta)]\|x_0-p_*\|}{\delta^{n+1}\prod_{k=0}^{n}[1-\alpha_k(1-\delta)]\|u_0-p_*\|} \\ &= \lim_{n\to\infty}\delta^{n+1}\frac{\prod_{k=0}^{n}[1-\alpha_k(1-\delta)][1-\beta_k(1-\delta)]}{\delta^{n+1}\prod_{k=0}^{n}[1-\alpha_k(1-\delta)]}\end{aligned}$$

Since we know that $x_0=u_0$ and $[1-\alpha_k(1-\delta)][1-\beta_k(1-\delta)]\le[1-\alpha_k(1-\delta)]$ , then

$$\lim_{n\to\infty}\frac{\|x_{n+1}-p_*\|}{\|u_{n+1}-p_*\|}\le\lim_{n\to\infty}\delta^{n+1}=0.$$

Thus,

$$\lim_{n\to\infty}\frac{\|x_{n+1}-p_*\|}{\|u_{n+1}-p_*\|}=0,$$

Then we conclude that $\{x_n\}_{n=0}^{\infty}$ *converges faster than* $\{u_n\}_{n=0}^{\infty}$.
Now, let us give a numerical example for analytic proof.

**Example 2.6** *Let us give the function* $T:[0,10]\to[0,10]$ *such that* $T(x)=\sqrt{x^2-8x+40}$. *It is easy to show that* $T$ *is an almost contraction on* $[0,10]$ *with fixed point* $p_*=5$. *Choose* $\alpha_n=\beta_n=\gamma_n=\frac{3}{4}$ *with the initial value* $x_0=1,99$.

*The following table shows that the new iteration method (0.6) converges faster than all Picard-Mann (0.5), Noor (0.4) and S (0.3) iteration methods.*

| $x_n$ | **The new iteration** | **Picard-Mann** | **Noor** | **S** |
|---|---|---|---|---|
| $x_1$ | 1,99000000000000 | 1,99000000000000 | 1,99000000000000 | 1,99000000000000 |
| $x_2$ | 4,99389167686581 | 4,92137488609813 | 4,17445437828765 | 5,01485375465129 |
| $x_3$ | 4,99996100255133 | 4,99389167686581 | 4,76149300996949 | 5,00164414119613 |
| $\vdots$ | $\vdots$ | $\vdots$ | $\vdots$ | $\vdots$ |
| $x_8$ | 5,00000000000000 | 4,99999998003362 | 4,99945966464220 | 5,00000002649966 |
| $\vdots$ | $\vdots$ | $\vdots$ | $\vdots$ | $\vdots$ |
| $x_{15}$ | 5,00000000000000 | 5,00000000000000 | 4,99999989115093 | 5,00000000000001 |
| $x_{16}$ | 5,00000000000000 | 5,00000000000000 | 4,99999996772625 | 5,00000000000000 |
| $\vdots$ | $\vdots$ | $\vdots$ | $\vdots$ | $\vdots$ |
| $x_{29}$ | 5,00000000000000 | 5,00000000000000 | 5,00000000000000 | 5,00000000000000 |

We notice that the new iteration method converges faster than the others, since it converges at the $8^{th}$ step while the Picard-Mann, Noor and *S* iteration methods converge at the $15^{th}$ step, $16^{th}$ step and $29^{th}$ step, respectively.
The following figure is a graphical presentation of the above results.

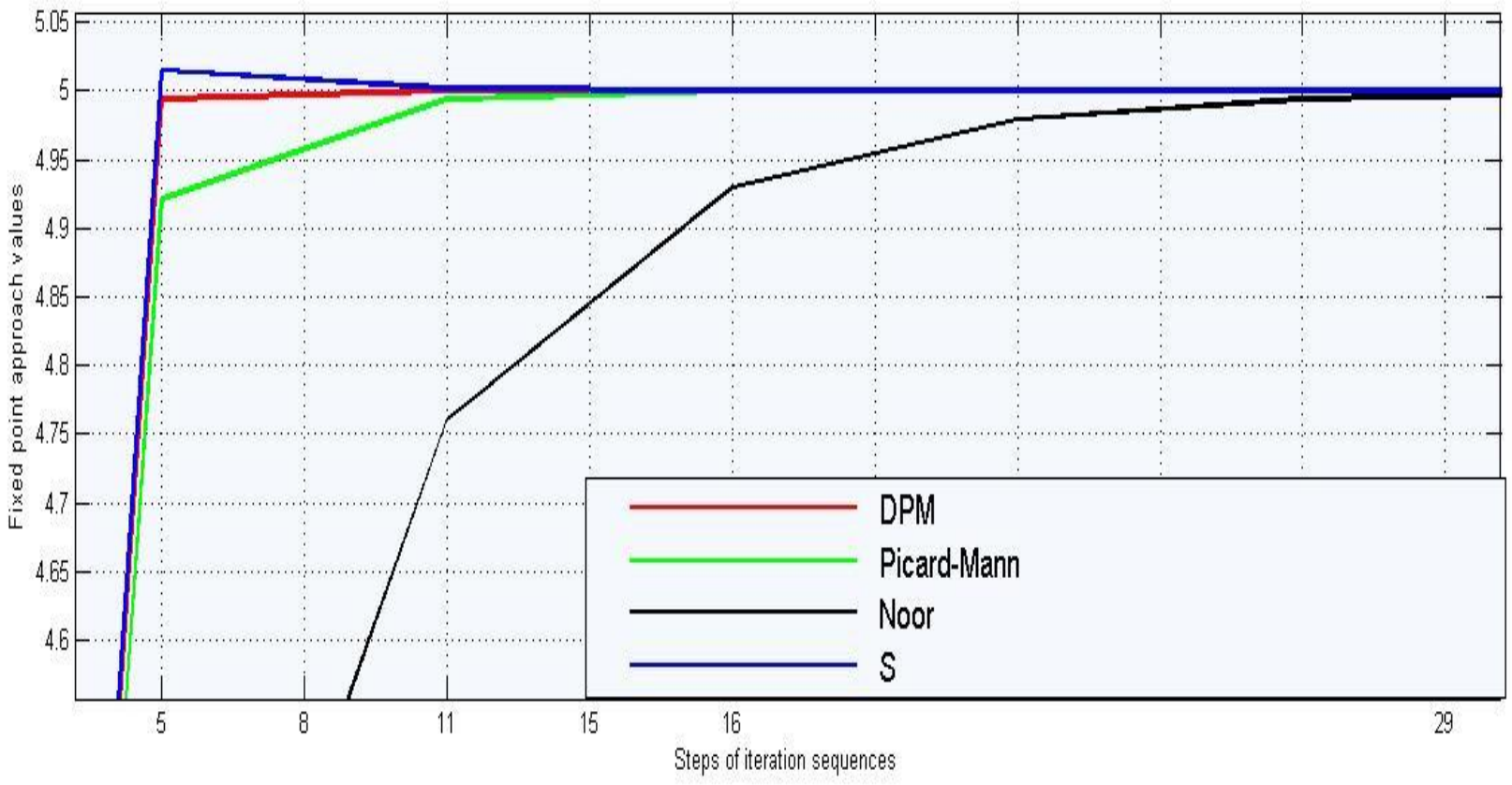

**Theorem 2.7** *Let* $S$ *be an approximate operator of* $T$ *. Let* $\{x_n\}_{n=0}^{\infty}$ *be an iterative sequence generated by* (0.6) *for* $T$ *and define an iterative sequence* $\{u_n\}_{n=0}^{\infty}$ *as follows:*

$$\begin{cases} u_0 \in C \\ u_{n+1} = S[(1-\alpha_n)v_n + \alpha_n S v_n] \\ v_n = S[(1-\beta_n)u_n + \beta_n S u_n],\ n \in \mathbb{N} \end{cases} \tag{1.24}$$

where $\{\alpha_n\}_{n=0}^{\infty}$ , $\{\beta_n\}_{n=0}^{\infty}$ be real sequences in $[0,1]$ satisfying (i) $\frac{1}{2} \le \alpha_n$ for all $n \in \mathbb{N}$ , If $Tp_* = p_*$ and $Sx_* = x_*$ such that $u_n \to x_*$ as $n \to \infty$ , then we have

$$\| p_* - x_* \| \le \frac{7\varepsilon}{1-\delta},$$

where $\varepsilon > 0$ is a fixed number.

**Proof.** From (0.6), (0.7) and (0.24), we have

$$
\begin{aligned}
\|y_n - v_n\| &\le \|T[(1-\beta_n)x_n + \beta_n Tx_n] - T[(1-\beta_n)u_n + \beta_n Su_n]\| \\
&\quad + \|T[(1-\beta_n)u_n + \beta_n Su_n] - S[(1-\beta_n)u_n + \beta_n Su_n]\| \\
&\le \delta(1-\beta_n)\|x_n - u_n\| + \delta\beta_n \|Tx_n - Su_n\| \\
&\quad + L(1-\beta_n)\|x_n - T[(1-\beta_n)x_n + \beta_n Tx_n]\| \\
&\quad + L\beta_n \|Tx_n - T[(1-\beta_n)x_n + \beta_n Tx_n]\| + \varepsilon \\
&\le \delta(1-\beta_n)\|x_n - u_n\| + \delta\beta_n \|Tx_n - Tu_n\| + \delta\beta_n \|Tu_n - Su_n\| \\
&\quad + L(1-\beta_n)\|x_n - Tx_n\| \\
&\quad + L(1-\beta_n)\|Tx_n - T[(1-\beta_n)x_n + \beta_n Tx_n]\| \\
&\quad + L\beta_n \|Tx_n - T[(1-\beta_n)x_n + \beta_n Tx_n]\| + \varepsilon \\
&\le \delta(1-\beta_n)\|x_n - u_n\| + \delta^2\beta_n \|x_n - u_n\| + \delta\beta_n L\|x_n - Tx_n\| + \delta\beta_n\varepsilon \\
&\quad + L(1-\beta_n)\|x_n - Tx_n\| \\
&\quad + L\delta\|x_n - (1-\beta_n)x_n - \beta_n Tx_n\| + L^2\|x_n - Tx_n\| + \varepsilon \\
&= \delta[1-\beta_n(1-\delta)]\|x_n - u_n\| \\
&\quad + L\{[1-\beta_n(1-\delta)] + L + \beta_n\delta\}\|x_n - Tx_n\| + \beta_n\delta\varepsilon + \varepsilon,
\end{aligned}
\tag{1.25}
$$

and

$$
\begin{aligned}
\|x_{n+1}-u_{n+1}\| &\le \|T[(1-\alpha_n)y_n+\alpha_n Ty_n]-T[(1-\alpha_n)v_n+\alpha_n Sv_n]\| \\
&\quad +\|T[(1-\alpha_n)v_n+\alpha_n Sv_n]-S[(1-\alpha_n)v_n+\alpha_n Sv_n]\| \\
&\le \delta\|(1-\alpha_n)(y_n-v_n)+\alpha_n(Ty_n-Sv_n)\| \\
&\quad +L\|(1-\alpha_n)y_n+\alpha_n Ty_n)-(1-\alpha_n+\alpha_n)T[(1-\alpha_n)y_n+\alpha_n Ty_n)]\| \\
&\quad +\varepsilon \\
&\le \delta(1-\alpha_n)\|y_n-v_n\|+\delta\alpha_n\|Ty_n-Tv_n\|+\delta\alpha_n\|Tv_n-Sv_n\| \\
&\quad +L(1-\alpha_n)\|y_n-T[(1-\alpha_n)y_n+\alpha_n Ty_n)]\| \\
&\quad +L\alpha_n\|Ty_n-T[(1-\alpha_n)y_n+\alpha_n Ty_n)]\|+\varepsilon \\
&\le \delta(1-\alpha_n)\|y_n-v_n\|+\delta^2\alpha_n\|y_n-v_n\|+\delta\alpha_n L\|y_n-Ty_n\|+\delta\alpha_n\varepsilon \\
&\quad +L(1-\alpha_n)\|y_n-Ty_n\| \\
&\quad +L(1-\alpha_n)\|Ty_n-T[(1-\alpha_n)y_n+\alpha_n Ty_n)]\| \\
&\quad +L\alpha_n\|Ty_n-T[(1-\alpha_n)y_n+\alpha_n Ty_n)]\|+\varepsilon \\
&\le \delta(1-\alpha_n)\|y_n-v_n\|+\delta^2\alpha_n\|y_n-v_n\|+\delta\alpha_n L\|y_n-Ty_n\|+\delta\alpha_n\varepsilon \\
&\quad +L(1-\alpha_n)\|y_n-Ty_n\|+L\delta\|y_n-(1-\alpha_n)y_n-\alpha_n Ty_n)\| \\
&\quad +L^2\|y_n-Ty_n\|+\varepsilon \\
&= \delta[1-\alpha_n(1-\delta)]\|y_n-v_n\| \\
&\quad +L\{[1-\alpha_n(1-\delta)]+L+\alpha_n\delta\}\|y_n-Ty_n\|+\alpha_n\delta\varepsilon+\varepsilon, \qquad (1.26)
\end{aligned}
$$

Substituting (0.25) in (0.26), we obtain

$$
\begin{aligned}
\|x_{n+1}-u_{n+1}\| &\le \delta^2[1-\alpha_n(1-\delta)][1-\beta_n(1-\delta)]\|x_n-u_n\| \\
&\quad +\delta[1-\alpha_n(1-\delta)]L\{[1-\beta_n(1-\delta)]+L+\beta_n\delta\}\|x_n-Tx_n\| \\
&\quad +\delta^2[1-\alpha_n(1-\delta)]\beta_n\varepsilon+\delta[1-\alpha_n(1-\delta)]\varepsilon \\
&\quad +L\{[1-\alpha_n(1-\delta)]+L+\alpha_n\delta\}\|y_n-Ty_n\|+\alpha_n\delta\varepsilon+\varepsilon
\end{aligned} \qquad (1.27)
$$

Since $\delta\in(0,1)$ and $\alpha_n,\beta_n\in[0,1]$ for all $n\in\mathbb{N}$, we have

$1-\alpha_n(1-\delta)<1,$

$1-\beta_n(1-\delta)<1,$

$\delta[1-\alpha_n(1-\delta)]<1,$

and using assumption we obtain

$1-\alpha_n\le\alpha_n.$

Hence, from (0.26), (0.27) and these inequalities, we have

$$\begin{aligned}
\|x_{n+1}-u_{n+1}\| &\le [1-\alpha_n(1-\delta)]\|x_n-u_n\| \\
&\quad +[1-\alpha_n+\alpha_n\delta]L\{[1-\beta_n(1-\delta)]+L+\beta_n\delta\}\|x_n-Tx_n\| \\
&\quad +(1-\alpha_n+\alpha_n)L\{[1-\alpha_n(1-\delta)]+L+\alpha_n\delta\}\|y_n-Ty_n\| \\
&\quad +(1-\alpha_n+\alpha_n)3\varepsilon+\alpha_n\varepsilon \\
&\le [1-\alpha_n(1-\delta)]\|x_n-u_n\| \\
&\quad +\alpha_n(1+\delta)L\{[1-\beta_n(1-\delta)]+L+\beta_n\delta\}\|x_n-Tx_n\| \\
&\quad +2\alpha_n L\{[1-\alpha_n(1-\delta)]+L+\alpha_n\delta\}\|y_n-Ty_n\|+7\alpha_n\varepsilon \\
&= [1-\alpha_n(1-\delta)]\|x_n-u_n\| \\
&\quad +\alpha_n(1-\delta)\left\{\frac{L(1+\delta)\{[1-\beta_n(1-\delta)]+L+\beta_n\delta\}\|x_n-Tx_n\|}{(1-\delta)}\right. \\
&\quad \left.+\frac{2L\{[1-\alpha_n(1-\delta)]+L+\alpha_n\delta\}\|y_n-Ty_n\|+7\varepsilon}{(1-\delta)}\right\}
\end{aligned}$$

Denote that

$$\begin{aligned}
a_n &= \|x_n-u_n\|, \\
\mu_n &= \alpha_n(1-\delta)\in(0,1) \\
\eta_n &= \frac{1}{(1-\delta)}\left\{\begin{matrix} L(1+\delta)\{[1-\beta_n(1-\delta)]+L+\beta_n\delta\}\|x_n-Tx_n\| \\ +2L\{[1-\alpha_n(1-\delta)]+L+\alpha_n\delta\}\|y_n-Ty_n\|+7\varepsilon \end{matrix}\right\}
\end{aligned}$$

Hence, from Lemma 1. 2, we have

$$\begin{aligned}
0 &\le \limsup_{n\to\infty}\|x_n-u_n\| \\
&\le \limsup_{n\to\infty}\left\{\frac{L(1+\delta)\{[1-\beta_n(1-\delta)]+L+\beta_n\delta\}\|x_n-Tx_n\|}{(1-\delta)}\right. \\
&\quad \left.+\frac{2L\{[1-\alpha_n(1-\delta)]+L+\alpha_n\delta\}\|y_n-Ty_n\|+7\varepsilon}{(1-\delta)}\right\} \\
&= \frac{7\varepsilon}{(1-\delta)}
\end{aligned}$$

We know from Theorem 2.1 that $x_n \to p_*$ and using hypothesis, we obtain $\|p_*-x_*\|\le\frac{7\varepsilon}{1-\delta}$.

## REFERENCES

**Agarwal, R., O Regan, D., Sahu, D.,** (2007), Iterative construction of fixed points of nearly asymptotically nonexpansive mappings. Journal of Nonlinear and Convex Analysis 8:61- 79.

**Berinde, V., (2007),** Iterative Approximation of Fixed Points, Springer, Berlin.

**Chugh, R., Kumar,V., Kumar, S.,** (2012), Strong Convergence of a New Three Step Iterative Scheme in Banach Spaces, *American Journal of Computational Mathematics*, Vol. 2 No. 4, pp. 345-357.

**Gürsoy, F., Karakaya, V., Rhoades, B. E.,** (2013), Data dependence results of new multi-step and S-iterative schemes for contractive-like operators, *Fixed Point Theory and Applications* Vol. 2013, (2013), doi:10.1186/1687-1812-2013-76.

**Ishikawa, S.,** (1974), Fixed Point By a New Iteration Method, *Proceedings of the American Mathematical Society*, Vol. 44, No.1, pp. 147-150.

**Karahan, I., Özdemir, M.,** (2013), A general iterative method for approximation of fixed points and their applications, *Advances in Fixed Point Theory*, Vol. 3, No.3, pp. 510-526.

**Karakaya, V., Doğan, K., Gürsoy, F., Ertürk, M.,** (2013), Fixed Point of a New Three-Step Iteration Algorithm under Contractive-Like Operators over Normed Spaces, *Abstract and Applied Analysis*, vol. 2013, Article ID 560258, 9 pages.

**Khan, S. H., (**2013), A Picard-Mann hybrid iterative process, Fixed Point Theory Appl. Article ID 69 (2013), doi:10.1186/1687-1812-2013-69.

**Mann, W.R.,** (1953), Mean Value Methods in Iteration, *Proceedings of the American Mathematical Society*, Vol. 4, No. 3, pp. 506-510.

**Osilike, M. O.,** (1995), Stability results for the Ishikawa fixed point iteration procedure. *Indian J. Pure Appl. Math*., 26 pp. 937-945.

**Pheungrattana, W.,.Suantai, R.,** (2011), On the Rate of convergence of Mann, Ishikawa, Nour and SP iterations for continuous on an Arbitrary interval, *Journal of computational and Applied Mathematics*, Vol. 235, No. 9pp. 3006-3914.

**Picard, E.** (1890), Memoire sur la theorie des equations aux derivees partielles et la methode des approximations successives". *J. Math. Pures Appl.* Vol. 6, No. 4, pp. 145-210.

**Soltuz, S.M., Grosan, T.,** (2008), Data dependence for Ishikawa iteration when dealing with contractive like operators", *Fixed Point Theory and Applications*, Vol. 2008**,** 7 pages.

**Weng, X.,** (1991), Fixed point iteration for local strictly pseudocontractive mapping, *Proc. Amer. Math. Soc.*, Vol. 113, pp. 727.